\documentclass[12pt]{amsart}
\usepackage{xypic}
\usepackage{epsfig}
\usepackage{xspace}
\usepackage{layout}
\usepackage{url}
\input xypic
\input xy
\xyoption{all}

\newcommand{\ie}{\textit{i.e.\xspace}}
\newcommand{\cform}[3]{\begin{array}{c}
{\scriptstyle #3}\\
#1\\
{\scriptstyle #2}\end{array}}

\def\qed{\hfill\square\medskip}
\newcommand{\beg}[2]{\begin{equation}\label{#1}#2\end{equation}}
\newcommand{\rref}[1]{(\ref{#1})}
\def\Proof{\medskip\noindent{\bf Proof: }}

\begin{document}

\newtheorem{thm}{Theorem}[section]
\newtheorem{conj}{Conjecture}[section]
\newtheorem{lem}[thm]{Lemma}
\newtheorem{cor}[thm]{Corollary}
\newtheorem{prop}[thm]{Proposition}
\newtheorem{rem}[thm]{Remark}

\theoremstyle{definition}
\newtheorem{defn}[thm]{Definition}
\newtheorem{examp}[thm]{Example}
\newtheorem{rmk}[thm]{Remark}
\def\r{\rightarrow}

\def\square{\hfill${\vcenter{\vbox{\hrule height.4pt \hbox{\vrule
width.4pt height7pt \kern7pt \vrule width.4pt} \hrule
height.4pt}}}$}

\newenvironment{pf}{{\it Proof:}\quad}{\square \vskip 12pt}
\date{\today}

\title[Laplaza sets]{Laplaza sets,
or how to select coherence diagrams for pseudo algebras}
\author{Thomas M. Fiore, Po Hu, and Igor Kriz}
\address{Thomas M. Fiore,
Department of Mathematics, University of Chicago, 5734 S.
University, Chicago, IL 60637, USA \\
and \\
Departament de Matem\`{a}tiques, Universitat Aut\`{o}noma de
Barcelona, 08193 Bellaterra (Barcelona), Spain}
\email{fiore@math.uchicago.edu}
\address{Po Hu, Department of Mathematics, Wayne State University, 1150
Faculty/Administration Building, 656 W. Kirby Detroit, MI 48202,
USA} \email{po@math.wayne.edu}
\address{Igor Kriz, Department of Mathematics, University of Michigan,
 2074 East
Hall, 530 Church Street, Ann Arbor, MI 48109-1043, USA}
\email{ikriz@umich.edu}
\thanks{
Thomas M.~Fiore was supported by NSF grant DMS 0501208 at the
University of Chicago. At the Universitat Aut\`{o}noma de Barcelona
he was supported by grant SB2006-0085 of the Spanish Ministerio de
Educaci\'{o}n y Ciencia under the Programa Nacional de ayudas para
la movilidad de profesores de universidad e investigadores
espa$\tilde{\text{n}}$oles y extranjeros. Po Hu was supported in
part by NSF grant DMS 0503814
and Igor Kriz was supported in part by NSF grant DMS 0305583. \\
\indent The authors thank Michael Shulman for some comments on an
ealier draft.}

\begin{abstract}
We define a general concept of pseudo algebras over theories and
$2$-theories. A more restrictive such notion was introduced in \cite{hu}, but as
noticed by M. Gould, did not capture the desired examples. The
approach taken in this paper corrects the mistake by introducing a
more general concept, allowing more flexibility in selecting
coherence diagrams for pseudo algebras.
\end{abstract}

\maketitle

\section{Introduction}

Generalizing algebras to pseudo algebras is a basic idea which has
recently become important in axiomatization of conformal field
theory \cite{fiore1,hu,hu1,lambda}, as well as in other subjects,
e.g.~\cite{blackwell}. While the exact settings vary, the kind of
algebras we are using generally have a set of {\em operations} and a
set of {\em
 identities}
(equations) the operations are required to satisfy. The
corresponding notion of {\em pseudo algebra} is a category rather
than a set. The operations are replaced by functors, and the
identities are replaced by natural isomorphisms which we call {\em
coherence isomorphisms}. Generally speaking, however, we now want
additional conditions, namely commutative diagrams which are to be
satisfied by the coherence isomorphisms. Such diagrams are generally
known as {\em coherence diagrams}. The question of what coherence
diagrams one should select is trickier than it may appear, and is
the main subject of this note.

\vspace{3mm} In \cite{hu}, the following scheme was suggested for
selecting coherence diagrams: take all diagrams which can be
``reasonably expected'' to commute. This means, take any word in our
algebra which can be formed by using variables represented by formal
symbols, and repeated use of operations which apply to them. Now
taking such a word, we can use the identities among the operations
(and any substitutions) to turn the word successively into other
words. It may happen, however, that one word $a$ can be turned into
another word $b$ in two different ways, using a different sequence
of identities. To such a situation, there corresponds in an obvious
way a coherence diagram (see for example \cite{maclane}). It was
suggested in \cite{hu} that {\em all} such diagrams should be
required to commute in a pseudo algebra.

\vspace{3mm} It turns out, however, that such a requirement is
unreasonably strong. For example, if the algebras in question are
commutative monoids where we denote the operation by $\oplus$, then
the word $a\oplus a$ can be turned into the same word either by an
empty sequence of identities,
or by an application (using substitution) of the identity
\beg{eop1}{a\oplus b =b\oplus a.} The corresponding coherence
diagram would then require that, in a pseudo commutative monoid, the
coherence isomorphism
$$\xymatrix{\tau_{ab}:a\oplus b \ar[r] & b\oplus a}$$
corresponding to \rref{eop1} satisfy \beg{eop2}{\tau_{aa}=Id.} This
however is unreasonably strong; we would like pseudo commutative
monoids to be the same thing as symmetric monoidal categories, and
those will not in general be equivalent to categories satisfying
\rref{eop2} (see Proposition \ref{pbad} below). The authors thank M.
Gould for this example, see Section 6 of \cite{gould2}.

\vspace{3mm} To correct this, one must generalize the notion of
pseudo algebra in a way that allows us to limit the scope of
coherence diagrams required, so that ``bad diagrams'' such as the
one mentioned above can be excluded. Surprisingly perhaps, as will
be shown in examples given below, there is not a single way to do
this which would cover all the examples desired. However, there is a
fairly simple and general scheme which includes all the cases needed
in \cite{hu, hu1, lambda}. This scheme amounts basically to
including coherence diagrams coming from processing one word $a$ to
another word $b$ using identities in the algebra, but with the
restriction that the formal variables used in each of the words $a$
and $b$ occur exactly once within each word in each identity. This
scheme requires an important restriction, namely all the identities
in the algebra must be between words which use each variable exactly
once. A precise formulation of this for the simplest case of
universal algebras modelled on one set (``$1$-sorted algebras'')
involves the language of operads, and its interplay with the
language of theories. The relevant concepts are defined in the next
section. The foundational results of \cite{fiore1} remain valid and
can be generalized to the new context, and hence all the substantive
results of \cite{hu, hu1, lambda} remain in effect. This idea is
also basically due to M. Gould. See \cite{gould1} and
\cite{gould2}.\footnote{The reference
 \cite{gould2} was
posted after the initial submission of the present article.}

\vspace{3mm} One kind of algebra which is of interest in conformal
field theory
 however is
the algebra of ``worldsheets'', i.e. Riemann surfaces with
analytically parametrized boundary components, and the operations of
disjoint union and gluing of boundary components of opposite
orientations. Such
 worldsheets
do not form a $1$-sorted algebra. This is because gluing requires
``dynamically indexed'' operations, in the sense that the possible
gluings depend on the set of boundary
 components
of the worldsheet and their orientations. Such structures are not
axiomatized by theories but by {\em $2$-theories}, introduced in
 \cite{hu}.
To apply the operad scheme for generating coherence diagrams in this
case, one needs to define {\em $2$-operads} (which we will do in
section \ref{s3} below).

\vspace{3mm} Even this operadic approach, however, is not
sufficiently general, since one is interested in algebras whose
identities do involve words with repeated symbols. Commutative
semi-rings give one such
 example,
where the distributivity axiom involves a word with recurring
variables on one side. In this case of pseudo commutative
semi-rings, the pseudo algebras should be symmetric bimonoidal
categories. The correct condition limiting coherence diagrams was
discovered by Laplaza \cite{laplaza}. The condition is somewhat
technical and will be explained in Example \ref{exlap} below. It is
not obvious what general scheme would select coherence diagrams
``correctly'' in accordance with what one expects for specific
examples of algebraic structures known. However, it is not difficult
to axiomatize what general {\em formal properties} such sets of
diagrams must satisfy. Such sets of diagrams, inspired by Laplaza's
diagrams, we call {\em Laplaza sets} in recognition of his
 contribution.

\section{Pseudo Algebras with Laplaza Sets in Theories}

Let us recall here the notion of a theory, which was first defined
in \cite{lawvere}. We will stick to the ``universal algebra'' point
of view, which is more advantageous for defining pseudo algebras.

\begin{defn}
Let $\Gamma$ be the category with objects $0,1,2,\dots$ where
$0=\emptyset$ and $k:=\{1, \dots , k\}$ for $k\geq1$. The morphisms
are maps, not necessarily order-preserving. Let
$\xymatrix@1{+:\Gamma \times \Gamma \ar[r] & \Gamma}$ be the
 functor
defined by $k+\ell:=\{1, \dots, k+\ell \}$ and by placing maps side
by side.
\end{defn}

\begin{defn} \label{theorydefinition}
A {\it theory} is a functor $\xymatrix@1{T:\Gamma \ar[r] & Sets}$
equipped with {\it compositions}
$$\xymatrix@1{\gamma: T(k) \times T(n_1) \times \dots \times T(n_k)
\ar[r] & T(n_1 + \dots + n_k)}$$ and a {\it unit} $1 \in T(1)$ such
that the following hold.
\begin{enumerate}
\item
\label{i1} The $\gamma$'s are {\it associative}, \ie
$$\gamma(w,\gamma(w^1,w^1_1,\dots,w^1_{n_1}),\gamma(w^2,w^2_1,\dots,w^2_{n_2}),
\dots, \gamma(w^k,w^k_1,\dots,w^k_{n_k}))=$$ $$\gamma(\gamma(w, w^1,
\dots, w^k), w^1_1, \dots,w_{n_1}^1,w^2_1, \dots, w^2_{n_2}, \dots,
w^k_1, \dots,w_{n_k}^k).$$
\item
\label{i2} The $\gamma$'s are {\it unital}, \ie $$\gamma(w,1,
\dots,1)=w=\gamma(1,w)$$ for all $w \in T(k)$.
\item
\label{i3} The $\gamma$'s are {\it equivariant} with respect to the
{\it functoriality} $()_f:=T(f)$ in the sense that
$$\gamma(w_f,w_1, \dots, w_{\ell})=\gamma(w,w_{f1}, \dots,
w_{fk})_{\bar{f}}$$ for every
function $\xymatrix@1{f:\{1, \dots, k\} \ar[r] & \{1, \dots, \ell
\}}$ where
$$\xymatrix{\bar{f}:\{1,2, \dots, n_{f1} + n_{f2} + \dots + n_{fk}
\} \ar[r] & \{1,2, \dots, n_1 + n_2 + \dots +
n_{\ell} \}}$$ is the
function that moves entire blocks according to $f$.
\item
\label{i4} The $\gamma$'s are {\it equivariant} with respect to
functoriality also in the sense that
$$\gamma(w,(w_1)_{g_1}, \dots, (w_k)_{g_k})= \gamma(w,w_1, \dots,
w_k)_{g_1 + \dots + g_k}$$ for all
functions $\xymatrix@1{g_i:\{1, \dots, n_i \} \ar[r] & \{1, \dots,
n_i'\}}$ where \\ $g_1 + \dots + g_k:\{1,2, \dots, n_1+ \dots + n_k \}
\rightarrow \{1,2, \dots, n_1 '+ \dots
+ n_k' \}$ is the
function obtained by placing $g_1, \dots, g_k$ next to each other
from left to right.
\end{enumerate}
The elements of $T(k)$ are called {\it words}.
\end{defn}

%\begin{examp}
If $X$ is a set, then $End(X)(n):=Map(X^n,X)$ defines the {\it
endomorphism theory} of $X$. Composition is the composition of
functions. The unit is $1_X$. If $\xymatrix@1{f:k \ar[r] & \ell}$
 is
a function and $w \in End(X)(k)$, then $w_f \in End(X)(\ell)$ is
defined by $w_f(x_1, \dots, x_{\ell})=w(x_{f1}, \dots, x_{fk})$.
This example allows us to define {\it algebras over theories}. A set
$X$ is a {\it $T$-algebra} when it is equipped with a morphism
$\xymatrix@1{T \ar[r] & End(X)}$ of theories.
%\end{examp}

\vspace{3mm} \noindent {\bf Remark:} Lawvere \cite{lawvere}
originally defined theories more elegantly as categories with the
set of natural numbers as objects, with the property that $k$ is a
categorical product of $k$ copies of $1$, with given projections.
The way this relates to Definition \ref{theorydefinition} is that
given a Lawvere theory $\mathcal{T}$, we define
$T(k)=\mathcal{T}(k,1)$. The axioms are then obviously satisfied. On
the other hand, given a theory $T$ in the sense of Definition
\ref{theorydefinition}, we set
\beg{elt}{\mathcal{T}(k,\ell)=T(k)^{\times\ell}.} (On the right hand
side, we mean the cartesian product of sets.) Composition of
elements of $\mathcal{T}(m,k)$ and $\mathcal{T}(k,\ell)$ is defined
by applying $\gamma$, which will get us to a product of $\ell$
copies of $T(mk)$, and then functoriality with respect to the map
$f:mk\r m$ which satisfies $f(i)\equiv i\mod m$.
%The $i$'th projection $k\r 1$ is the injection $1\mapsto i\in k$
%composed with the unit in $T(1)$, and the unit is the product
%of the $i$'th projections, $i=1,...,k$.
To obtain the $i$'th projection $k\r 1$ in $\mathcal{T}$, we
substitute the injection $1\mapsto i\in k$ into the unit $1$ in
$T(1)$. The unit in the category $\mathcal{T}$ on the object $k$ is
the product of the $i$'th projections, $i=1,...,k$. One needs the
equivariance axioms (\ref{i3}), (\ref{i4}) to prove associativity and
unitality, although one can show that the axioms have some
redundancy (i.e. can be deduced from special cases). (For complete
detail, see Chapter 6 of \cite{fiore1}.)

Theories model {\em $1$-sorted universal algebras}. By this we mean
algebras whose definition calls for one set with some operations
required to satisfy certain prescribed identities. Then the set of
words $T(n)$ is the set of all operations in $n$-symbols that arise
as a composite of finitely many basic operations in $T$. The symbols
in words are allowed to repeat.

More generally, an $n$-sorted (or $I$-sorted, where $I$ is an
indexing
 set)
algebra calls for $I$ sets and operations which are allowed to apply
to prescribed sets in $I$, and produce an element of another
prescribed set of $I$. Again, prescribed identities (equations) are
required to hold. For example, a ring and a module form a $2$-sorted
algebra.

\vspace{3mm} An important observation about categories $\mathcal{C}$
of all multi-sorted algebras with given operations and identities is
that if we have the category of all multi-sorted algebras
 $\mathcal{D}$
whose operations and identities form (possibly empty) subsets of the
sets of operations and identities of $\mathcal{C}$, then there is a
forgetful functor $\mathcal{C}\r\mathcal{D}$ which has a left
adjoint. We usually refer to the left adjoint as the functor $F$
taking the {\em free $\mathcal{C}$-algebra on a
$\mathcal{D}$-algebra}. To prove the existence of these left
adjoints, we note two constructions standard in algebra, the
construction of a {\em free algebra} and the construction of a {\em
quotient}. The first is the special case of left adjoint to the
forgetful functor to systems of sets. To construct a free
multi-sorted algebra of the given kind on a system of sets, take the set
of all formal words using the operations on the elements of the
applicable sets, and then factor out by the smallest equivalence
relation which includes all the required identities and is preserved
by operations (an operation on equivalent elements gives equivalent
results). The quotient construction, for a multi-sorted algebra $X$,
and a relation $\sim$ on its elements, gives a universal {\em
quotient} of $X$ which is a multi-sorted algebra with the same
operations and identities, and in which any pair of elements $x\sim
y$ are identified. Once again, it is constructed by taking the
smallest equivalence relation which contains the given relation and
is preserved by operations. To construct the left adjoint $F$
mentioned above, we let $F(X)$ be the free $\mathcal{C}$-algebra on
$X$, and take the quotient under the relation identifying all
$\mathcal{D}$-words in elements of $X$ with their result in $X$.

$I$-sorted algebras are not axiomatized by Lawvere theories,
although the formalism can be adapted to them. We do not take this
approach here, however, since we will need an even more general
context, described in the next section.

\begin{defn}
\label{operaddefinition} The notion of {\em operad}
is defined by following verbatim
Definition \ref{theorydefinition}, except that we replace the
morphisms of the category $\Gamma$ by all bijections,
and restrict in the equivariance axioms (\ref{i3}),
(\ref{i4}) to all bijective maps. Algebras over
an operad are defined precisely analogously as we defined algebras
over a theory.
\end{defn}

\vspace{3mm} One advantage of Definition \ref{theorydefinition} is
that it exhibits the fact that the notion of a theory is itself an
$\mathbb{N}$-sorted algebra where $\mathbb{N}$ is the set of all
natural numbers. By the above remarks, then, we have forgetful
functors from theories to operads, to
sequences of sets. We call a sequence of sets $Z=\{Z(n)\}_{n \geq0}$ a {\it collection}.
A theory is free on an operad
if and only if it is generated by operations and equations where the
equations involve no repetition of variables on either side and
the exact same variables occur on both sides: then the
underlying operad consists of all words in chosen
variables $a_1,\dots,a_n$ which can be written, using the operations,
where each variable has to be used exactly once.

These statements are proved easily. The point is, in both cases, we
can define the operad generated by
the given operations and identities. The free theory on those
operations modulo these equations is the same thing as the free
theory on the operad free on the
operations modulo these equations: maps in both directions are
exhibited and proved inverse to each other by the universal
properties.

\vspace{3mm} These definitions worked in the category of sets.
However, the category $Cat$ of small categories and functors has
properties analogous to those of the category of sets and maps. One
therefore immediately gets the analogous notion of {\it internal
theory} in $Cat$. We call such an internal theory a {\it categorical
theory}. Explicitly, to define a categorical theory, replace in
Definition \ref{theorydefinition} $Sets$ by $Cat$; all the axioms
(1)-(4) can be rewritten as diagrams, which are the same in $Cat$ as
in $Sets$. In particular, if $X$ is a category, then we can
similarly define a functor $n \mapsto End_{Cat}(X)(n)=Funct(X^n,X)$
from $\Gamma$ to $Cat$ (the morphisms in $End_{Cat}(X)(n)$ are
natural transformations). If $T$ is a categorical theory, then a
category $X$ is a {\it $T$-algebra} when it is equipped with a
morphism $\xymatrix@1{T \ar[r] & End_{Cat}(X)}$ of categorical
theories.

\vspace{3mm} It is useful to also note that categorical theories $T$
have an alternate description: both $Obj(T)$ and $Mor(T)$ are
theories, while source, target, and identity are morphisms of
theories. This makes $Mor(T)\times_{Obj(T)}Mor(T)$ into a theory,
and composition is a morphism of theories. Finally, let us note that
we may consider the category of graphical pre-theories. A {\it
graphical pre-theory} $T$ is a theory $\{Obj(T)(n)\}_{n \geq0}$ and a
collection $\{Mor(T)(n)\}_{n \geq0}$ together with maps $Source$,
$Target$, and $Id$ which satisfy the usual unital property, but
without composition. Equivalently, a graphical pre-theory
$(Obj(T),Mor(T))$ is an internal reflexive graph in the category of
collections with the additional
property that the object collection $Obj(T)$ is a theory. From
another point of view, however (reinterpreting graphs and categories
with a fixed object set as multi-sorted algebras over $Sets$),
categorical theories and graphical pre-theories are also
multi-sorted algebras in $Sets$. We then have, using our general
observations about multi-sorted algebras, a forgetful functor from
the category of categorical theories to the category of graphical
pre-theories. This functor has a left adjoint, which is the free
categorical theory on a graphical pre-theory. Further, both of these
functors preserve the object collection.

\vspace{3mm} With the notions of graphical pre-theory and categorical
theory in hand, we are now ready to introduce pseudo algebras over a
theory with respect to a Laplaza set. The purpose of a Laplaza set
$S$ is to specify which diagrams are required to commute in a pseudo
$T$-algebra, and to do this we force certain diagrams to commute in
an associated categorical theory $T^{\prime}_{S}$. Roughly speaking,
a Laplaza set for a theory $T$ is a set $S$ of words in $T$, and a
diagram of coherence isomorphisms is required to commute whenever
the words of the source and target are in the Laplaza set $S$. To
make this precise, we construct from an ordinary theory $T$ a
categorical theory $T^{\prime}_{S}$ which has an isomorphism between
each free composite of words in $T$ and their composite in $T$. Some
diagrams of isomorphisms will commute, exactly which ones is decided
by the Laplaza set $S$. Any morphism $\xymatrix@1{T^{\prime}_{S}
\ar[r] & End_{Cat}(X)}$ of categorical theories then takes the
abstract isomorphisms and commutative diagrams to coherence
isomorphisms and coherence diagrams for $X$, thus, such a morphism
is a pseudo $T$-algebra with respect to the Laplaza set $S$.

\vspace{3mm}
\begin{defn}
\label{dl1} A {\em Laplaza set} for a theory $T$ is an arbitrary
collection of sets $S(n)\subseteq T(n)$. For a given Laplaza set,
define a categorical theory $T^{\prime}_{S}$ as follows. First
define a graphical pre-theory $G_S$. The theory $Obj(G_S)$ is the
free theory on the collection $\{S(n)\}_{n\geq0}$. The set
$Mor(G_S)(n)$ contains one arrow $\iota_{a,b}$ between each pair of
objects in $Obj(G_S)(n)$ which map to the same element of
$S(n)\subseteq T(n)$, and $Id_a=\iota_{a,a}$. Now $T^{\prime}_{S}$
is the quotient of the free categorical theory $F_S$ on the
graphical pre-theory $G_S$ by the relations
\beg{el1}{\iota_{ab}\sim\iota_{ba}^{-1},}
\beg{el2}{\parbox{3.5in}{If $\alpha,\beta\in Mor(F_S)(n)$ satisfy
$Source(\alpha)=Source(\beta)=x$, $Target(\alpha)=Target(\beta)=y$
and $x,y$ project to the same element of $S(n)\subseteq T(n)$, then
$\alpha=\beta$.}} A {\em pseudo $(T,S)$-algebra} (or any other
permutation of these words, e.g. a pseudo algebra over $T$ with
respect to the Laplaza set $S$ etc.) is a $T^{\prime}_{S}$-algebra,
i.e. a morphism $\xymatrix@1{T^{\prime}_{S} \ar[r] & End_{Cat}(X)}$
of categorical theories.
\end{defn}

The special case considered in \cite{hu,fiore1} is $S=T$. The
 difficulty
discovered by M. Gould is expressed very strongly by the following

\begin{prop}
\label{pbad} Let $T$ be the theory of commutative monoids. Then
every pseudo $(T,T)$-algebra $A$ is equivalent to a strictly
symmetric monoidal category, i.e. a
 category
$A^{\prime}$ which is a strict algebra over the theory of
commutative monoids.
\end{prop}

\Proof Select representatives $a_i$, $i\in I$, of isomorphism
classes of $A$, and assume $I$ is a linearly ordered set, with
minimum $0$. Assume $a_0$ is the unit. Let the operation be
$\oplus$. Then define the category $A^\prime$, with operation $+$,
to have objects $a_i$ with $i\in I$ as well as formal sums \beg{ela}{a_{i_1}+
a_{i_2}+\cdots+a_{i_n}} where $0<i_1\leq\cdots \leq i_n$. The operation $+$
on elements of the form \rref{ela} with $n$ and $m$ summands,
respectively, is the sum of the form \rref{ela} which shuffles the
elements together so that the indices are again in non-decreasing
order. The sum $x+a_0$ is defined to be $x$ for any $x$. Clearly,
this operation is commutative, associative, and
 unital.

To define morphisms, define a map $F$ from the proposed objects of
 $A^\prime$
to $Obj(A)$ by sending \rref{ela} to
$$a_{i_1}\oplus...\oplus a_{i_n}$$
and $a_0$ to $a_0$. Then pull back $Mor(A)$ via this map, thus
 promoting
$F$ into an equivalence of categories. It remains to define the
 operation
$+$ on morphisms. To define $f+g$, consider $F(f)\oplus F(g)$ and
compose on both sides with any coherence isomorphisms needed to
shuffle the source and target back to order; then pull back via $F$.
The result is unique by the
 observation
that the switch coherence iso $a\oplus a\rightarrow a\oplus a$ must
be the identity. \qed

All the formal results of \cite{fiore1} generalize to pseudo
algebras with
Laplaza sets, in particular these structures form $2$-categories in
the obvious way, and enjoy pseudo limits and bicolimits. The proofs
of \cite{fiore1} work essentially word by word. From this point of
view, pseudo algebras over $T$ in the sense of \cite{fiore1} are
pseudo algebras with respect to the Laplaza set $S=T$, i.e. a
special case. For the biadjunctions discussed in \cite{fiore1}, the
appropriate forgetful functor is associated with a morphism of
theories
$$\xymatrix{\phi: T_1\ar[r] & T_2}$$
together with Laplaza sets $S_i\subseteq T_i$, satisfying the
condition
$$\phi(S_1)\subseteq S_2.$$
Then there is a forgetful $2$-functor from pseudo
$(T_2,S_2)$-algebras to pseudo $(T_1,S_1)$-algebras which enjoys a
left biadjoint constructed by the same method as in \cite{fiore1}.

\vspace{3mm} If $T$ is the free theory on an operad
$C$, then there is a canonical example of a
Laplaza set associated with $C$
which is often useful, namely the collection $\{C(n)\}_{n\geq0}$ itself.
One notes that if we denote by $\Gamma$ (resp. $\Sigma$) the category
 whose objects
are natural numbers and morphisms $m\r n$ are maps (resp. bijections)
$\{1,...,m\}\r \{1,...,n\}$, then the free theory $T$
on an operad $C$ is
\beg{eggg}{\Gamma\times_\Sigma C=\cform{\coprod}{n\geq 0}{}
\Gamma(n,?)\times_{\Sigma_n}C(n),
}
so the canonical map $C(n)\r T(n)$ is injective.

If $C$ is the operad defining commutative monoids and $T$ is the
free theory on this operad, then we do have an inclusion
$C\subset T$, and pseudo $(T,C)$-algebras are
unbiased symmetric monoidal categories. An {\it unbiased} symmetric
monoidal category is one in which
$n$-fold products $x_1 \otimes x_2 \otimes \cdots  \otimes x_n$ are chosen
in addition to binary products, and there are accompanying coherence isomorphisms satisfying the obvious
coherence diagrams. The Laplaza set $C$ essentially by definition forces precisely
the coherence diagrams of \cite{maclane} and their unbiased counterparts. The equivalence of the category of
unbiased symmetric monoidal categories with the category of symmetric monoidal categories is
essentially the coherence theorem of Mac Lane.

\vspace{3mm} Let $C$ be an operad. Then define a {\it categorical operad}
$C^\prime$ (i.e. operad internal in categories)  as follows: the
objects are the free operad on the collection $\{C(n)\}_{n\geq0}$,
and there is precisely one isomorphism between any two objects of
$C^\prime$ which map to the same element of $C$. A {\it pseudo
$C$-algebra} is a $C'$-algebra, \ie a morphism $\xymatrix@1{C'\ar[r]
& End_{Cat}(X)}$ of categorical operads. This is the way that
 pseudo
algebras over operads are defined in \cite{gould1} and
\cite{gould2}.

\begin{prop}
\label{pgood} Let $T$ be a theory which is free on an operad $C$,
and let $S$ be the corresponding Laplaza set.
Then a pseudo $(T,S)$-algebra is the same thing as a pseudo
$C$-algebra.
\end{prop}

\Proof Let $T^{\sharp}_{S}$ be the free categorical theory on the
categorical operad $C^\prime$. There is an obvious morphism of
categorical theories \beg{epgo2}{\xymatrix{T^{\sharp}_{S}\ar[r] &
T^{\prime}_{S}.} } Indeed, this map is obtained by universality and
the observation that the relations in Definition \ref{dl1} imply the
operad relations in $C^\prime$. We claim that \rref{epgo2} is an
isomorphism. The object theory of $T^{\sharp}_{S}$ is the free
theory on the free operad $Obj(C')$ on the collection
$\{C(n)\}_{n\geq0}$. This is the same as the object theory of $T^{\prime}_{S}$,
namely the free theory on the collection $\{C(n)\}_{n\geq0}$. It is
easy to see from the definitions that \rref{epgo2} must be full. To
show that it is faithful, we claim that we can use the universality of $T^{\prime}_{S}$ to
construct a left inverse. This is equivalent to the following
statement: \beg{eggs}{\parbox{3.5in}{In $T^{\sharp}_{S}$, there is
precisely one isomorphism between any two objects which under the
canonical map
 to
$T$ project to the same element of $C$.}
}
To prove \rref{eggs}, the objects of $T^{\sharp}_{S}$ are of the form
$(f,u)$ where $f$ is a function, and $u$ is an element of the free
 operad
on the collection $\{C(n)\}_{n\geq0}$. Denote by $\overline{u}$ the image of $u$
 in
$T$. Note first of all that if $\overline{u}=\overline{v}$, there is
 certainly
an isomorphism between $(f,u)$ and $(f,v)$ in $T^{\sharp}_{S}$.
Similarly, if $\sigma$ is a bijection, then $(f\sigma,u)$ is equal to
$(f,u_\sigma)$ in $T^{\sharp}_{S}$, so there is an identity isomorphism
between them. Therefore, by induction one sees that when
$\overline{u}_f=\overline{v}_g$ in $T$,
then $(f,u)$ is isomorphic to $(g,v)$ in $T^{\sharp}_{S}$.

Therefore, the only non-trivial statement in \rref{eggs} is uniqueness,
and it suffices to assume that the two objects concerned are the same
 object
of $C^\prime$.
In other words, we must prove that for $u\in Obj C^\prime$, the only
 self-map
$u\r u$ in $T^{\sharp}_{S}$ is the identity. Let, therefore, $(f,a)$ be
 another
such self-map where $f$ is a function and \beg{eggs0}{\xymatrix{a:v\ar[r] & w}} is a
morphism in $C^\prime$. Then we must have in particular
\beg{eggs1}{(f,v)=(f,w)=u} in $Obj(T^{\sharp}_{S})$. But
$Obj(T^{\sharp}_{S})$ is a free theory on a collection.
Therefore, the elements of $Obj(T^{\sharp}_{S})$ are formal words we
 can
write in a given ordered set of variables using the words in
 $\{C(n)\}_{n \geq 0}$. Repetition
of variables is allowed, and there is no identification. The word will
 belong
to $Obj C^\prime$ (the free operad) if no repetition of variables
 occurs
and each variable is used exactly once. Thus, $f$ must be a
bijection, and hence we may as well assume $f=Id$. It then follows
that $v=w$ and hence $a$ is the identity. \qed

\begin{examp}
\label{exlap} The previous proposition is not sufficient to define
symmetric bimonoidal categories. In this case, the theory $T$ is the
theory of commutative semi-rings. This is not a free theory on an
operad (the reason being that
distributivity involves repetition of symbols on one side of the
equation). This is the case \cite{laplaza} where the original Laplaza set was
defined: one lets $S(n)$ consist of all words which, when converted
to the form of a sum of monomials using distributivity, identifying
a monomial $m$ with $1\cdot m$ and deleting any $0$ summands
($0\cdot m=0$), reduce to a sum of distinct square free monomials
(monomials which are permutations of each other are considered
equal). With this choice of Laplaza set $S$, a pseudo $(T,S)$-algebra is an
unbiased symmetric bimonoidal category. The equivalence of the category
of unbiased symmetric bimonoidal categories with the category of
symmetric bimonoidal categories is essentially Laplaza's coherence
theorem in \cite{laplaza}.
\end{examp}

\section{Laplaza Sets in $2$-theories}

\label{s3}

The notion of 2-theory is defined in \cite{hu}. The main example of
interest here is the $2$-theory of commutative monoids with
cancellation. We recapitulate these definitions here before turning
to Laplaza sets for $2$-theories.

\begin{defn}
\label{d2t} A {\it 2-theory} consists of a natural number $k$, a
theory $T$, sets
$$\Theta(w;w_1, \dots, w_n)$$
for all $w_1, \dots, w_n, w \in T(m)^k, m\geq0$, and the following
{\it operations}.
\begin{enumerate}
\item\label{operation1}
For each $w \in T(m)^k$ there exists a {\it unit} $1_w \in
End(X)(w;w)$.
\item\label{operation2}
For all $w, w_i, w_{ij} \in T(m)^k$ there is a function called {\it
$\Theta$-composition}.

\begingroup
\vspace{-2\abovedisplayskip} \small
$$\gamma:\Theta(w;w_1, \dots, w_q) \times \Theta(w_1;w_{11}, \dots,
w_{1p_1}) \times \cdots \times \Theta(w_q;w_{q1}, \dots, w_{qp_q})$$
$$\xymatrix{\ar[r] & \Theta(w;w_{11}, \dots, w_{qp_q})}$$
\endgroup
\noindent
\item\label{operation3}
Let $w, w_1, \dots, w_q \in T(m)^k$. For any function $\iota:\{1,
\dots, p\} \rightarrow \{1, \dots, q\}$ there is a function
$$\xymatrix@1{()^{\iota}:\Theta(w;w_{\iota(1)}, \dots,w_{\iota(p)})
\ar[r] & \Theta(w;w_1, \dots, w_q)}$$ called {\em
 $\Theta$-functoriality}.
\item\label{operation4}
Let $w, w_1, \dots, w_q \in T(m)^k$. For any function $f:\{1,
\dots,m\} \rightarrow \{1, \dots, \ell \}$ there is a function
$$\xymatrix@1{()_f:\Theta(w;w_1, \dots, w_q) \ar[r] & \Theta(w_{f};
(w_1)_{f}, \dots, (w_q)_{f})}$$ where $w_{f}$ means to substitute
$f$ in each of the words in the $k$-tuple $w$. This function is
called {\it $T$-functoriality}.
\item\label{operation5}
For $u_i \in T(k_i), i=1, \dots, m$ and $w, w_1, \dots, w_q \in
T(m)^k$ let $v_j:=\gamma^{\times k}(w_j;u_1^{\times k},
\dots,u_m^{\times k})$ for $j=1, \dots, q$ and furthermore let
$v:=\gamma^{\times k}(w;u_1^{\times k}, \dots,u_m^{\times k})$. Then
there is a function
$$\xymatrix@1{(u_1,\dots,u_m)^*:\Theta(w;w_1, \dots, w_q) \ar[r] &
\Theta(v;v_1, \dots, v_q)}$$ called {\it $T$-substitution}. Here
$\gamma^{\times k}$ means to use the composition of the theory $T$
in each of the $k$ components, which coincides with composition in
the theory $T^k$ with $T^k(m):=T(m)^k$.
\end{enumerate}
These operations satisfy the following {\it relations} (cf. pages
152-154 of \cite{fiore1}):
\begin{enumerate}
\item\label{ax1}
$\Theta$-composition is associative and unital in an analogous sense
as \rref{i1} and \rref{i2} in the definition of a theory
\item\label{ax2}
$\Theta$-functoriality is functorial in the sense that for functions
$$\diagram
\{1,...,p\}\rto^\iota &\{1,...,q\}\rto^\theta &\{1,...,r\},
\enddiagram$$
we have $()^\theta()^\iota=()^{\theta\iota}$ and $()^{Id}=Id$.
\item\label{ax3}
$\Theta$-composition is equivariant with respect to
$\Theta$-functoriality in two ways, analogously as \rref{i3} and
\rref{i4} in the definition of a
 theory.
\item\label{ax4}
\label{axt} $T$-functoriality is functorial in the sense that for
functions
$$\diagram
\{1,...,n\}\rto^f &\{1,...,m\}\rto^g &\{1,...,\ell\},
\enddiagram$$
we have $()_g ()_f=()_{gf}$ and $()_{Id}=Id$.
\item\label{ax5}
$T$-substitution is compatible with composition and unit in the
sense that if $w,w_1,...,w_q\in T(m)^k$, $t_i\in T(k_i) $ $s_{ij}\in
T(k_{ij}) $, $1\leq i \leq m$, $1\leq j\leq k_i$, if we set
$r_i=\gamma^{\times k} (t_{i}^{\times k}, s_{i1}^{\times
k},...s_{ik_i}^{\times k}) $,
$$(r_1,...,r_m)^*=(s_{11},...,s_{mk_m})^*(t_1,...,t_m)^*,$$
and also
$$(1,...,1)^*=Id.
$$
\item\label{ax6}
$\Theta$-composition is $T$-equivariant in the sense that if
$f:\{1,...,m\}\r\{1,...,\ell\}$ is a function, $w,w_i,w_{ij}\in
T(m)^k$, $\alpha\in\Theta(w;w_1,...,w_q) $,
$\alpha_j\in\Theta(w_j;w_{j1},...,w_{jp_j})$ for $j=1,...,q$, then
$$\gamma(\alpha_f;(\alpha_1)_f,...,(\alpha_q)_f)=
\gamma(\alpha;\alpha_1,...,\alpha_q)_f.$$
\item\label{ax7}
$\Theta$-functoriality and $T$-functoriality commute in the sense
that for functions $\iota:\{1,...,p\}\r\{1,...,q\}$ and
$f:\{1,...,m\}\r\{1,...\ell\}$, we have
$(\alpha^\iota)_f=(\alpha_f)^\iota$ for all $\alpha \in
\Theta(w;w_{\iota(1)}, \dots, w_{\iota(p)})$.
\item\label{ax8}
$\Theta$-functoriality and $T$-substitution commute:
$$(u_1,...,u_m)^*()^\iota=()^\iota(u_1,...,u_m)^*.$$
\item\label{ax9}
$T$-functoriality and $T$-substitution commute in the sense that for
$u_i\in T(k_i)$, and $f_i:\{1,...,k_i\}\r\{1,...,k_{i}^{\prime}\}$,
if we denote by $f:\{1,...,\sum k_i\}\r\{1,...,\sum
k_{i}^{\prime}\}$ the juxtaposition of the functions $f_i$, then
$$()_f(u_1,...,u_m)^*=((u_1)_{f_1},...,(u_m)_{f_m})^*.
$$
\item\label{ax10}
$T$-functoriality and $T$-substitution also commute in the sense
that if $f:\{1,...,m\}\r\{1,...\ell\}$ is a function, then
$$(u_1,...,u_\ell)^*()_f=()_{\bar{f}}(u_{f(1)},...,u_{f(m)})^*
$$
\item\label{ax11}
\label{axn} $T$-substitution and $T$-composition commute in the
sense that if $u_i\in T(k_i)$, $i=1,...,m$, $\alpha\in
\Theta(w;w_1,...,w_q)$, $\alpha_\ell\in \Theta(w_\ell;w_{\ell
1},...,w_{\ell p_\ell})$ for $\ell=1,...,q$, and
$\beta=(u_1,...,u_m)^*\alpha$, $\beta_i=(u_1,...,u_m)^*\alpha_i$,
then
$$(u_1,...,u_m)^*\gamma(\alpha;\alpha_1,...,\alpha_q)=
\gamma(\beta;\beta_1,...,\beta_q).$$
\end{enumerate}
\end{defn}

\vspace{3mm} Again, there is an alternate ``Lawvere-style''
categorical description (cf. \cite{hu}). More precisely, in that
sense, a $2$-theory consists of a natural number $k$, a theory $T$,
and a (strict) contravariant functor $\Theta$ from $T$ to the
category of small categories (and functors) with the following
properties. Let $T^k$ denote the category with the same objects as
$T$ (natural numbers) and such that
$Hom_{T^k}(m,n)=(Hom_{T}(m,n))^{\times k}$. Then for every natural
number $m$,
$$Obj(\Theta(m))=\cform{\coprod}{n\geq0}{}Hom_{T^{k}}(m,n),$$
for every morphism $\phi:m\r n$ in $T$ the map $Obj(\Theta(n))\r
Obj(\Theta(m))$, which is a part of $\Theta(\phi)$, is given by
precomposition with $(\phi,...,\phi)$, and lastly every
$$\psi\in Hom_{T^{k}}(m,n)$$
is the product, in $\Theta(m)$, of the $n$-tuple
$$w_{1},...,w_{n}\in Hom_{T^{k}}(m,1)$$
with which it is identified by the fact that $T$ is a theory.

This ``categorical'' definition is shown to be equivalent with the
``algebraic'' Definition \ref{d2t} as follows. Operations
\rref{operation1}, \rref{operation2}, \rref{operation3} and
relations \rref{ax1}, \rref{ax2}, \rref{ax3} are equivalent to
saying that $\Theta(m)$ is a category, similarly as in the case of
theories. The key point is to identify
$$Hom_{\Theta(m)}(\prod_{i=1}^n w_i,w)=
\Theta(w;w_1,...,w_n).$$ Operations \rref{operation4} and
\rref{operation5} are the morphism part of the strict 2-functor
$\Theta:T^{op} \r Cat$. Relations \rref{ax4}, \rref{ax5},
\rref{ax9}, and \rref{ax10} are then equivalent to saying that the
morphism part of $\Theta$ is a functor into sets. Relations
\rref{ax6}, \rref{ax7}, \rref{ax8}, \rref{ax11} are then equivalent
to promoting $\Theta$ to a functor into the category of small
categories and functors.

\vspace{3mm} Roughly speaking, the point of $2$-theories is to index
algebras with ``dynamically indexed'' operations. We have an algebra
$I$ over a certain theory and sets $X_{(i_1,...,i_k)}$ where
$i_1,...,i_k\in I$. The kind of $n$-ary operations we allow on the
$X$'s take as input tuples of elements
$$x_{j}\in X_{(w_{j1},...,w_{jk})}$$
and produce an element of
$$X_{(w_1,...,w_k)}$$
where $w_{ji}$, $w_{i}$ are certain specified words, all in the same
 given
set of variables. Relations can also be specified on these
operations, leading to the above definition.

\vspace{3mm} More formally, if $I$ is a set and $\xymatrix@1{X:I^k
\ar[r] & Sets}$ is a map, then we have a 2-theory $End(X)$ where
$End(I)$ is the underlying theory, and $End(X)(w;w_1, \dots, w_n)$
is the set of maps
$$\xymatrix@1{X \circ w_1 \circ d^m \times \cdots \times X \circ w_n
\circ d^m \ar@{->}[r] & X \circ w \circ d^m}$$ for maps
$\xymatrix@1{w_1, \dots, w_n, w:(I^m)^k \ar[r] & I^k}$. Here
$\xymatrix@1{d^m:I^m \ar[r] & (I^m)^k}$ is the diagonal map. Using
this example, a map $\xymatrix@1{X:I^k \ar[r] & Sets}$ is a {\it
$(\Theta,T)$-algebra} when it is equipped with a morphism
$$\xymatrix@1{(\Theta,T) \ar[r] & (End(X),End(I))}$$ of 2-theories.

\vspace{3mm} It is useful to note again that the notions in the last
paragraph are defined in the category of sets and maps, but can be
defined in $Cat$ when we replace ``sets'' by ``categories'' and
``maps'' by ``functors''. In particular, associated to any category
$I$
 and any
strict 2-functor $\xymatrix@1{X:I^2 \ar[r] & Cat}$ there is a {\it
categorical $2$-theory} \beg{eend}{(End_{Cat}(X),End_{Cat}(I))} and
an {\it algebra over a categorical 2-theory} is defined as a
morphism from a given categorical 2-theory to \rref{eend}. It is
also useful to note that, again, we can alternately define
categorical $2$-theories as consisting of an object $2$-theory and
morphism $2$-theory, which satisfy the axioms of a category, but in
the category of $2$-theories.

\vspace{3mm} $2$-theories are not multi-sorted algebras in the usual
sense. However, if we have already fixed a theory $T$, then
$2$-theories over $T$ are multi-sorted algebras (sorted over
$(n+1)k$-tuples of elements of $T(m)$ for all $m$). Therefore, we
can speak of a free $2$-theory on a system of sets
$\Xi(\gamma;\gamma_1,...,\gamma_n)$ over a given theory $T$. We can
also, once $T$ and $\Xi$ are fixed, impose equivalence relations
$\sim$ on each of the sets $\Xi(\gamma;\gamma_1,...,\gamma_n)$.
There exists a universal quotient of $\Xi$ which forms a $2$-theory
over $T$ and on which $\sim$ will turn into equality.

\vspace{3mm}
Given a {\it categorical theory} $T$, we may define the notion of
{\it graphical pre-$2$-theory} $(\Xi,T)$. This consists of a $2$-theory
$(Obj(\Xi),Obj(T))$ with underlying theory $Obj(T)$, a set
$Mor(\Xi)(\gamma;\gamma_1,...,\gamma_n)$ for all
 $\gamma;\gamma_1,...,\gamma_n$
$k$-tuples of words  of the theory $Mor(T)$, as well as $Source$,
 $Target$, $Id$ maps
satisfying the usual unitality axioms, but no composition. There are no
 2-theory axioms
on the sets $Mor(\Xi)(\gamma;\gamma_1,...,\gamma_n)$. The free
categorical $2$-theory on the graphical pre-$2$-theory $(\Xi,T)$ has
 the same
underlying categorical theory $T$ and the same object $2$-theory
$(Obj(\Xi),Obj(T))$.

%\vspace{3mm} Considering a categorical theory $T$ we may define a
%graphical $2$-theory which specifies $2$-theories $Obj(\Xi)$ over
%$Obj(T)$ and $Mor(\Xi)$ over $Mor(T)$ with $Source$, $Target$, $Id$
%satisfying the usual unitality axioms, but no composition. We can
%then again speak of a free categorical $2$-theory on this data.
%Similarly as above, we will define a graphical pre-$2$-theory.
%We will require that there be an underlying graphical theory, and
%that the object set be a $2$-theory on the objects of this graphical
%theory. Finally, we will have sets $\Xi(\gamma;\gamma_1,...,\gamma_n)$
%where $\gamma_i$, $\gamma$ are morphisms in the underlying graphical
%theory, but no $2$-theory axioms on these sets.

\begin{defn}
A {\em Laplaza set} $(\Sigma,S)$ for a $2$-theory $(\Xi,T)$ consists
of an arbitrary collection of sets $S(n)\subseteq T(n)$ and an
arbitrary system of sets
$\Sigma(\gamma;\gamma_1,...,\gamma_n)\subseteq
\Xi(\gamma;\gamma_1,...,\gamma_n)$ where
$\gamma;\gamma_1,...,\gamma_n$ are words in $S(m)^k$ for some $m$.
For a given Laplaza set, define a categorical $2$-theory
$(\Xi^{\prime}_{\Sigma},T^{\prime}_{S})$ as follows. First,
$T^{\prime}_{S}$ was already defined in Definition \ref{dl1}. Next,
define $Obj(\Xi^{\prime}_{\Sigma})$ as the free $2$-theory on the
system of sets $\Theta(\delta;\delta_1,...,\delta_n)$ indexed by
words $\delta;\delta_1,...,\delta_n$ in $Obj(T^{\prime}_{S})(m)^k$
that project to some $\gamma;\gamma_1,...,\gamma_n$ in $S(m)^k$.
Here the set $\Theta(\delta;\delta_1,...,\delta_n)$ is equal to
$\Sigma(\gamma;\gamma_1,...,\gamma_n)$ where
$\gamma;\gamma_1,...,\gamma_n$ is the projection of
$\delta;\delta_1,...,\delta_n$ to $S(m)^k$.

To define $Mor(\Xi^{\prime}_{\Sigma})$, first define (as in the
theory case) a graphical pre-$2$-theory $\Gamma_{(\Sigma,S)}$ over
the categorical theory $T_S^{\prime}$ whose objects are
$Obj(\Xi^{\prime}_{\Sigma})$, and morphisms are one $\iota_{a,b}$
between each $a,b\in Obj(\Xi^{\prime}_{\Sigma})$ indexed over tuples
$\delta;\delta_1,...,\delta_n$,
$\epsilon;\epsilon_1,...,\epsilon_n$, respectively, indexed over the
$\iota_{\delta_{i},\epsilon_{i}}$ for $i=\emptyset, 1,..,n$, {\em
when, additionally,} $a$, $b$ project to the same word in
$\Sigma(\gamma;\gamma_1,...,\gamma_n)
\subseteq\Xi(\gamma;\gamma_1,...,\gamma_n)$. Again, we impose that
$\iota_{a,a}=Id$ to get a reflexive graph.

Now we define $\Xi^{\prime}_{\Sigma}$ as the quotient of the free
categorical $2$-theory $\Phi_{(\Sigma,S)}$ on $\Gamma_{(\Sigma,S)}$
by the relations \beg{e2l1}{\iota_{ab}\sim\iota_{ba}^{-1},}
\beg{e2l2}{\parbox{3.5in}{If $\alpha,\beta\in
Mor(\Phi_{(\Sigma,S)})$ satisfy $Source(\alpha)=Source(\beta)=x$,
$Target(\alpha)=Target(\beta)=y$ and $x,y$ project to the same
element of $\Sigma$, then $\alpha=\beta$.}} (Note that since
elements of $\Sigma$ are only allowed to be indexed over tuples of
words in $S^k$, $\alpha$ and $\beta$ are necessarily indexed over the
same tuple of morphisms.) A {\em pseudo $(\Xi,T,\Sigma, S)$-algebra}
(or any other permutation of these words, e.g. a pseudo algebra over
$(\Xi, T)$ with respect to the Laplaza set $(\Sigma,S)$ etc.) is an
algebra over the categorical 2-theory
$(\Xi^{\prime}_{\Sigma},T^{\prime}_{S})$, i.e. a morphism
$\xymatrix@1{(\Xi^{\prime}_{\Sigma},T^{\prime}_{S}) \ar[r] &
 (End_{Cat}(X),End_{Cat}(I))}$ of categorical 2-theories.
\end{defn}

Again, pseudo algebras over a $2$-theory $(\Xi,T)$ with respect to a
given Laplaza set $(\Sigma,S)$ enjoy pseudo limits; the proofs of
\cite{fiore1} generalize easily. Therefore, we can speak of stacks
of pseudo $(\Xi,T,\Sigma,S)$-algebras.

\begin{defn}
One can define the notion of {\it $2$-operad} by repeating
Definition \ref{d2t} with the following changes.
\begin{itemize}
\item
$T$ is the free theory on an operad $C$.
\item
The word ``function'' is replaced by ``bijection'' in
$\Theta$-functoriality and also in axioms pertaining to
$\Theta$-functoriality as appropriate.
\end{itemize}
In particular, the indexing words of a 2-operad can be theory words.
\end{defn}

\vspace{3mm} By an analogous argument as before, the forgetful
functor from $2$-theories to $2$-operads has a left adjoint, the
free $2$-theory $(\Xi,T)$ on a 2-operad $(\Delta,C)$. In this case,
$(\Delta,C)$ provides a canonical choice of Laplaza set in
$(\Xi,T)$.

It is worth commenting again that a $2$-theory $(\Xi,T)$ is free
over a $2$-operad when $T$ is free over an operad
$C$, and $\Xi$ can be expressed as a quotient of a
free $2$-theory on a set of given generating operations in tuples of
words in $T$ by equations both sides of which
have the exact same indeterminates, which do not
repeat. The proof is analogous to the case of the
adjunction between theories and operads.

The main example of a 2-operad of interest here is the {\it 2-operad
of commutative monoids with cancellation}. In this example, $k=2$,
and $T$ is the theory of commutative monoids. We describe this
2-operad via its algebras.
\begin{defn}
\label{edefw} A strict 2-functor $\xymatrix@1{X:I^2 \ar[r] & Cat}$
is an {\it algebra over the 2-operad of commutative monoids with
cancellation} if $I$ is an algebra over the operad of commutative
monoids, and $X$ is equipped with natural functors
$$\xymatrix@1{+:X_{a,b} \times X_{c,d} \ar[r] & X_{a+c, b +d} }$$
$$\xymatrix@1{\check{?}:X_{a+c,b+c} \ar[r] & X_{a,b}}$$
$$0 \in X_{0,0}$$
satisfying the following axioms.
\begin{enumerate}
\item
The operation $+$ is {\it commutative}. $$
\xymatrix@R=3pc@C=3pc{X_{a,b} \times X_{c,d} \ar[r]^{+} \ar[d] &s in
X_{a+c, b +d} \ar@{=}[d]
\\ X_{c,d} \times X_{a,b} \ar[r]_{+} & X_{c+a,d+b}}$$
\item
The operation $+$ is {\it associative}.
$$\xymatrix@C=5pc@R=3pc{ (X_{a,b} \times X_{c,d}) \times X_{e,f} \ar[d]
\ar[r]^{+ \times 1_{X_{e,f}}} & X_{a+c,b+d} \times X_{e,f} \ar[d]^+
\\ X_{a,b} \times (X_{c,d} \times X_{e,f})
\ar[d]_{1_{X_{a,b}}\times +} & X_{(a+c)+e,(b+d)+f} \ar@{=}[d]
\\ X_{a,b} \times X_{c+e,d+f} \ar[r]_+ & X_{a+(c+e),b+(d+f)}
 }$$
\item
The operation $+$ has {\it unit} $0 \in X_{0,0}$.
$$\xymatrix@R=3pc@C=3pc{X_{a,b} \times \{0\} \ar[r]^{ +  }
\ar[dr]_{pr_1} & X_{a+0,b+0} \ar@{=}[d]
\\ & X_{a,b}} $$
\item
The operation $\check{?}$ is {\it transitive}.
$$\xymatrix@C=3pc@R=3pc{X_{(a+c)+d, (b+c)+d} \ar@{=}[d]
\ar[r]^-{\check{?}} & X_{a+c,b+c}
\ar[d]^{\check{?}} \\
X_{a+(c+d),b+(c+d)} \ar[r]_-{\check{?}} & X_{a,b} }$$
\item
The operation $\check{?}$ {\it distributes} over the operation $+$.
$$\xymatrix@R=3pc@C=3pc{X_{a+c,b+c} \times X_{e,f} \ar[r]^{+}
\ar[dd]_{\check{?} \times 1_{X_{e,f}}} & X_{(a+c)+e,(b+c)+f}
\ar@{=}[d] \\ & X_{(a+e)+c,(b+f)+c} \ar[d]^{\check{?}}
\\ X_{a,b} \times X_{e,f}
\ar[r]_{+} & X_{a+e,b+f} }$$
\item
Trivial cancellation is trivial.
$$\xymatrix@1@C=3pc@R=3pc{X_{a+0,b+0} \ar[r]^-{\check{?}} \ar@{=}[dr]
 &
X_{a,b} \ar[d]^{1_{X_{a,b}}} \\ & X_{a,b}}$$
\end{enumerate}
\end{defn}

\vspace{3mm} \noindent {\bf Remark:} The reader should note that
characterizing commutative monoids with cancellation in terms of its
algebras is purely a matter of language. In terms of Definition
\ref{d2t}, $T$ is the theory of commutative monoids, $k=2$, and the
$\Theta$'s are all the operations we can express by iterating the
operations $+$, $\check{?}$ and $0$ using iteration and substitution
without repetition of variables. Such operations are identified
subject to the relations (1)-(6).

\vspace{3mm} In this example,
$$+\in End(X)((pr_1+pr_3, pr_2+pr_4);(pr_1,pr_2),(pr_3, pr_4))$$
for $pr_i\in End(I)(4)$ and $$\check{?}\in
End(X)((pr_1,pr_2);(pr_1+pr_3,pr_2+pr_3))$$ for $pr_i \in
End(I)(3)$. Projection to the $i$-th coordinate is the same as
subsitution by the injective map $\xymatrix@1{\iota_i:\{1\} \ar[r] &
\{1,\dots, k\}}$, $\iota_i(1)=i$. This is the reason why theory words
are permitted as indexing words in the definition of 2-operad.

It is clear that the theory of commutative monoids is given by an
operad. The full force of the $T$-functoriality of $2$-operads as we
defined it is used in the transitivity axiom (4): In the composition
of cancellations \beg{eforce}{X_{(a+c)+d,(b+c)+d}\r X_{a+c,b+c}\r
X_{a,b},} the first map is obtained from the cancellation
$$X_{(u+d,v+d)}\r X_{u,v}$$
by $T$-substituting $a+e$ for $u$, $b+h$ for $v$, but then applying
$T$-functoriality with respect to the map identifying the variables
$e$ and $h$ into $c$.

\begin{examp}
Recall that a {\it worldsheet} is a real, compact, not necessarily
connected, two dimensional manifold with complex structure and
analytically parametrized boundary components.

\begin{prop}
\label{pworld} Worldsheets form a pseudo algebra over the 2-theory
of commutative monoids with cancellation (with Laplaza set
corresponding to the $2$-operad described above, on which the
$2$-theory is free).
\end{prop}

\Proof Let $I$ denote the symmetric monoidal category of finite sets
and bijections with $+=\coprod$. For finite sets $A$ and $B$,
$X_{A,B}$ is the category of worldsheets with inbound components
labelled by $A$ and outbound components labelled by $B$. A morphism
in $X_{A,B}$ is a holomorphic diffeomorphism that preserves boundary
parametrizations and boundary component labellings. If $f$ and $g$
are bijections the functor $X_{f,g}$ corresponds to boundary
relabellings. The operation $+$ is the disjoint union of
worldsheets, $\xymatrix@1{\check{?}:X_{a+c,b+c} \ar[r] & X_{a,b}}$
is the self gluing of boundary components with the same label in
$c$, and $0 \in X_{0,0}$ is the empty manifold. The coherence
isomorphisms
 from
the previous definition are defined by noting that we have canonical
embeddings $X,Y\r X\coprod Y$ and a canonical map $X\r \check{X}$
which is an embedding on the interior $X-\partial X$ of $X$. We see
then for $n$ distinct worldsheets $X_1,\dots,X_n$, and any
worldsheet $X$ obtained by repeated use of $+$ and $\check{?}$ on
$X_1,\dots,X_n$ where we use each $X_i$ exactly once, there are
canonical maps $X_i\r X$
 which
are embeddings on the interior of $X_i$ with disjoint images, whose
 union
is dense in $X$. Further, these embeddings commute with the
coherence isomorphisms corresponding to the identities in Definition
\ref{edefw}. Therefore, any coherence diagram corresponding to two
ways of
 processing
a word on distinct variables $X_1,\dots,X_n$ into another word using
 identities
in Definition \ref{edefw} will commute on the union of images of the
 $X_i$'s
in the result $X$ of the composite operation. But this union is
dense in $X$. \qed

More strongly, worldsheets actually form a {\em stack} of pseudo
algebras over the 2-operad of commutative monoids with cancellation:
the construction of the stack structure given in \cite{hu} is
correct in this new definition.

\end{examp}

\begin{defn}
A {\it conformal field theory} (in the most abstract sense) is a
morphism of stacks of pseudo algebras over the 2-operad of
commutative monoids with cancellation.
\end{defn}

%\bibliographystyle{hplain}
%\bibliography{biblio}

\def\cprime{$'$}

\end{document}